
\input amstex
\documentstyle{amsppt}
\refstyle{A}

\def\cC{{\Cal C}}
\def\cD{{\Cal D}}\def\cE{{\Cal E}}\def\cF{{\Cal F}}\def\cG{{\Cal G}}
\def\cH{{\Cal H}}\def\cL{{\Cal L}}\def\cP{{\Cal P}}

\def\bbR{{\Bbb R}}\def\bbC{{\Bbb C}}
\def\bbP{{\Bbb P}}
\def\w{{\mathchoice{\,{\scriptstyle\wedge}\,}{{\scriptstyle\wedge}}
{{\scriptscriptstyle\wedge}}{{\scriptscriptstyle\wedge}}}}

\def\euso{{\frak so}}
\def\SO{\text{SO}}
\def\ts{\textstyle }
\def\>{\medspace}
\def\ldot{\,.\,}

\topmatter
\title On extremals with prescribed Lagrangian densities
\endtitle
\author
Robert L. Bryant
\endauthor
\affil
Duke University
\endaffil
\address
Department of Mathematics,
Duke University,
PO Box 90320,
Durham, NC 27708-0320
\endaddress
\date
May 10, 1994
\enddate

\email
bryant\@math.duke.edu
\endemail

\keywords
Lagrangians, minimal surfaces, isosystolic metrics, harmonic maps
\endkeywords

\subjclass
53A10, 58E20
\endsubjclass

\dedicatory
Dedicated to Professor Eugenio Calabi on the occasion of his seventieth
birthday.
\enddedicatory

\thanks
This work was begun during a visit to Oberwohlfach in June 1993, continued
during the conference in honor of Professor Calabi in August 1993 at the Scuola
Normale Superiore in Pisa, and completed during a visit to MSRI in Spring~1994.
 The support of these institutions and that from the National Science
Foundation in the form of Grant DMS-9205222 is hereby gratefully acknowledged.
\endthanks

\abstract
Consider two manifolds~$M^m$ and $N^n$ and a first-order Lagrangian $L(u)$ for
mappings $u:M\to N$, i.e., $L$ is an expression involving $u$ and its first
derivatives whose value is an $m$-form (or more generally, an $m$-density)
on~$M$.  One is usually interested in describing the extrema of the functional
$\Cal L(u) = \int_M L(u)$, and these are characterized locally as the solutions
of the Euler-Lagrange equation~$E_L(u)=0$ associated to $L$.  In this note I
will discuss three problems which can be understood as trying to determine how
many solutions exist to the Euler-Lagrange equation which also satisfy $L(u) =
\Phi$, where $\Phi$ is a specified $m$-form or $m$-density on~$M$.  The first
problem, which is solved completely, is to determine when two minimal graphs
over a domain in the plane can induce the same area form without merely
differing by a vertical translation or reflection.  The second problem,
described more fully below, arose in Professor Calabi's study of extremal
isosystolic metrics on surfaces.  The third problem, also solved completely, is
to determine the (local) harmonic maps between spheres which have constant
energy density.
\endabstract

\endtopmatter

\document

\head 0. Introduction \endhead

In the study of {\smc pde}, the search for solutions with special properties
has a long history.  Adding extra conditions, such as symmetry assumptions,
frequently reduces the complexity of a problem to a point where it can be
explicitly solved.  The resulting special solutions are then often used as the
basis for formulating conjectures and as comparisons to the general solution so
as to gain information about global behavior, asymptotics, long time existence,
etc.  In another direction, one is often interested in uniqueness questions,
i.e., determining how much extra information must be specified in order to
determine a unique solution to the {\smc pde} under study.

In this note a combination of these two approaches is considered.  I am going
to describe three problems in which the {\smc pde} under consideration is the
Euler-Lagrange equation for some first order functional and the extra condition
takes the form of specifying the value of the Lagrangian density of the
solution in advance.  Normally this extra condition so overdetermines the
problem that there will be no solutions at all.  What I will be interested in
will be the special cases for which, not only do solutions exist, but
uniqueness fails.

Let me describe a simple, concrete case of the general problem.  Suppose given
a first order functional of the form
$$
\cL(z) = \int_\cD\ L(x,y,z,z_x,z_y)\, dx \w dy
\tag 1
$$
where $\cD\subset\bbR^2$ is a plane domain and $z$ is a (real-valued) function
on this domain.  For any given function~$F$ on~$\cD$, the problem is to
describe the space of solutions to the Euler-Lagrange equation for~$\cL$ which
also satisfy the first order condition
$$
L(x,y,z,z_x,z_y) = F(x,y).
\tag 2
$$
Some basic insight into the nature of this problem can be seen as follows.
Differentiating~(2) with respect to~$x$ and~$y$ yields the equations
$$
\aligned
L_x + L_z\,z_x + L_p\,z_{xx} + L_q\,z_{xy} &= F_x\\
L_y + L_z\,z_y + L_p\,z_{xy} + L_q\,z_{yy} &= F_y\\
\endaligned
\tag 3
$$
while the Euler-Lagrange equation of~(1) is of the form
$$
 L_{pp}\,z_{xx} + 2L_{pq}\,z_{xy} + L_{qq}\,z_{yy} = \tilde L
\tag 4
$$
where $\tilde L$ only involves $x$, $y$, $z$, $p=z_x$ and $q=z_y$.  The three
equations (3-4) can be solved for the second partials (i.e., the 2-jet) of~$z$
in terms of the lower order partials (i.e., the 1-jet) of~$z$ provided that the
obvious determinant does not vanish (which will, in general, be the case).  The
usual methods of total differential equations~\cite{BCG} then imply that there
is at most a three-parameter family of solutions to these three equations.
Coupled with the original first-order equation~(2), this results in a maximum
of two parameters worth of solutions to the system.

However, for the generic choice of~$F$, the system~(2-4) will be incompatible
and will have no solutions at all.  The question is to determine those~$F$ for
which this system is sufficiently compatible that there will exist more than
one solution.  The ``extra'' compatibility of this system will turn out to be
expressed as a system of partial differential equations on~$F$, which need not
have any solutions.

The first two problems to be treated below were suggested to me by Professor
Calabi and are precisely of the type described above.  In the first, the
Lagrangian is the area Lagrangian for a graph~$z = u(x,y)$ over a plane domain,
i.e.,
$$
L = \sqrt{1+p^2+q^2}.
$$
Thus, the problem is to determine when two minimal graphs can determine the
same area function.  Obviously, two minimal graphs which differ by a vertical
translation or a vertical reflection will induce the same area function, so all
of these solutions should be regarded as equivalent.  In \S1, I show that,
aside from the obvious case of two planes of the same slope, there are
essentially three other cases of pairs of inequivalent minimal graphs which
induce the same area form.  All of these turn out to be known minimal surfaces,
with analytic continuations to complete embedded minimal surfaces in space,
though the last family appears to have been considered only recently~\cite{Ka}.

The second problem is like the first, but for the Lagrangian
$$
L = \frac{2pq}{\sqrt{(2pq)^2-(p^2+q^2-1)^2}}\,,
$$
which arose in Professor Calabi's treatment of extremal isosystolic metrics on
surfaces~\cite{Ca2}.  Here, the answer is not as complete.  I have only been
able to show that, up to some simple equivalences, there is at most a three
parameter family of pairs of inequivalent solutions of the Euler-Lagrange
equation which induce the same Lagrangian.  Though the analysis could, in
principle, be completed without the introduction of new ideas, the calculation
which must then be carried out is beyond the range of any computer to which I
have access.  The motivation for this analysis is explained in the introduction
to Section~2.

The third and final problem to be considered here is of a different nature,
though it is in the same general class.  Given a (simply-connected piece of a)
Riemannian $n$-manifold~$\Sigma^n$ of constant sectional curvature~$K$, what
are the possible constant values that can be assumed by the energy
density~$||df||^2$ of a harmonic mapping $f:\Sigma\to S^N$?  I show that, if
$K$ is positive, and~$||df||^2$ is constant, then $||df||^2 = m(m{+}n{-}1)K$
for some integer~$m\ge 0$ and that, moreover,such an $f$ is, in a suitable
sense, the restriction to~$\Sigma$ of harmonic polynomials of degree~$m$
on~$\bbR^{n+1}$.  Thus, in this setting local harmonic mappings of constant
energy density extend to global harmonic mappings.  This result generalizes
results of Wallach~\cite{W} and myself~\cite{Br}.  Moreover, the non-uniqueness
of~$f$ in this case is reduced to a representation-theoretic problem analogous
to that considered in~\cite{DW} and~\cite{DZ} in the case of minimal isometric
immersions.

\head 1. Minimal graphs with equal areas \endhead

I want to begin by considering a simple, but nonetheless interesting example,
namely the Lagrangian for the Euclidean area of a graph
$(x,y)\mapsto\bigl(x,y,u(x,y)\bigr)$,
$$
L(x,y,z,p,q) = \sqrt{1+p^2 + q^2}.
$$
The Euler-Lagrange equation for the associated area functional is
$$
{{\partial\hfil}\over{\partial x}}
\left({{u_x}\over{\sqrt{1+u_x^2 + u_y^2}}}\right)
+
{{\partial\hfil}\over{\partial y}}
\left({{u_y}\over{\sqrt{1+u_x^2 + u_y^2}}}\right)
= 0.
\tag 1
$$

Given an open, connected domain~$\cD$ in the $xy$-plane, and a (sufficiently
differentiable) function~$F$ on~$\cD$ which satisfies $F\ge1$, the problem is
to determine how many functions~$u$ in~$\cD$ satisfy the minimal surface
equation as well as the first order {\smc pde}
$$
1 + u_x^2 + u_y^2 = F^2.
\tag 2
$$
If $u$ is one solution of~(1-2), then all of the functions of the form~$\pm u +
c$ (where $c$ is a constant) are also solutions of~(1-2), but I will regard
these as equivalent to the original solution.  The problem is to find those
functions~$F$ for which there exist at least two inequivalent solutions, i.e.,
two minimal surfaces graphed over the same domain~$\cD$ in the plane which
induce the same area function~$F$ on~$\cD$ {\it without\/} being congruent by
vertical translation or reflection.

Before stating the result, let me make a few comments about some evident
symmetries of the problem.  If $u$ is a solution of the minimal surface
equation on a domain~$\cD$, then for any constant $\lambda>0$, the function
$u_\lambda(x,y) = \lambda\,u(x/\lambda,y/\lambda)$, defined on the scaled
domain~$\lambda\cD$ is also a solution of the minimal surface equation (its
graph is the $\lambda$-dilation in space of the graph of~$u$).  Moreover if
$L(u) = F$, then $L(u_\lambda) = F_\lambda$ where $F_\lambda(x,y) =
F\bigl(x/\lambda,\,y/\lambda\bigr)$.  Thus, if $F$ allows two inequivalent
solutions, then so will $F_\lambda$.  I will say that the functions~$F_\lambda$
are equivalent up to scale.  A similar remark applies to translations or
rotations of the domain~$\cD$ in the $xy$-plane.  It obviously suffices to
solve the problem up to these sorts of trivial equivalences.

\proclaim{Theorem~1}
Suppose that a function~$F$ on a connected plane domain~$\cD$ is the area
function of two inequivalent minimal graphs over~$\cD$.  Then, after possibly
translating, scaling, and rotating in the $xy$-plane, $F$ is the restriction of
one of the following functions to a subdomain:
\roster
\item $F\equiv c$ where $c>1$ is a constant and~$\cD=\bbR^2$.  There is
1-parameter family of corresponding inequivalent minimal graphs, namely, the
circle of planes~$\Sigma_\psi$ defined by the equations
$$
z=\sqrt{c^2-1}\,\bigl(\cos\psi\,\,x+\sin\psi\,\,y\bigr).
$$
\item $F = \coth x$ and $\cD = \{ (x,y)\ \vrule\ x>0\}$.  There is a
1-parameter family of corresponding inequivalent minimal graphs, namely the
graphs~$\Sigma_\psi$ each of which is the part lying over~$\cD$ of Scherk's
fifth \rom(minimal\rom) surface
$$
\sinh x\,\sinh z = \cos(y + \psi).
$$
\item The function~$F_\phi$ on the
domain~$\cD_\phi=\{ (x,y)\ \vrule\ x^2+y^2 > \sin^2\phi\}$ defined by
$$
F_\phi(x,y) = \sqrt{{x^2+y^2+\cos^2\phi}\over{x^2+y^2-\sin^2\phi}} \>.
$$
There are exactly two inequivalent \rom(multi-valued\rom) minimal graphs
associated to~$F_\phi$.  They are the upper and lower sheets of the
``heli-catenoid'' described below.
\item The function~$F_{a,c}$ defined
on~$\cD_{a,c}=\{ (x,y)\ \vrule\ a\,\cosh x + c\,\cos y >1\}$ by
$$
F_{a,c}(x,y) = \sqrt{{a\,\cosh x+c\,\cos y+1}\over{a\,\cosh x+c\,\cos y-1}}\>,
$$
where the constants $a$ and $c$ satisfy $|a{-}c|<1<a{+}c$.  There are exactly
two inequivalent minimal graphs over~$\cD_{a,c}$ associated to~$F_{a,c}$.  Up
to vertical translation and reflection, they are subsets of a single complete,
doubly-periodic minimal surface in~$\bbR^3$, as described below.
\endroster
\endproclaim

\noindent The remainder of the section is devoted to the proof of this theorem.

\subhead 1.1. Reduction to an equation \endsubhead
First, a couple of preliminary remarks.  If a continuously differentiable
function~$u$ satisfies the minimal surface equation~(1) on a connected
domain~$\cD\subset\bbR^2$, then~$u$ is known to be real-analytic.  If $u$ also
satisfies~(2), then it follows that $F$ must be real analytic as well.  Thus,
it may assumed from the start that $F$ is real-analytic.

Moreover, I may also suppose that the locus~$F=1$ consists of isolated points,
since, otherwise, either $F\equiv1$ (in which case it is clear that the only
solutions to~(2) are $u\equiv c$) or else the set~$F=1$ would contain an
analytic arc~$\alpha\subset\cD$.  The graph of~$u$ would then be tangent to a
(horizontal) plane along~$\alpha$.  However, two connected minimal surfaces
which are tangent along a curve are identical, so the graph of $u$ would be a
horizontal plane, implying that $F\equiv1$, which, as has been seen, does not
furnish a solution to the problem.  Thus, by excising isolated points where
$F=1$, I may assume from now on that $F>1$ at every point of~$\cD$.

Set $F = \coth\mu$ where $\mu$ is a positive (analytic) function.%
\footnote{Just why the equations that follow are drastically simpler in terms
of~$\mu$ than they would be if expressed directly in terms of~$F$ is a mystery
to me.}
Then, for any solution~$u$ of~(1-2), there must exist a function~$\theta$
on~$\cD$ (well-defined up to an additive multiple of~$2\pi$) so that
$$
u_x = {{\cos\theta}\over{\sinh\mu}}
\qquad\text{and}\qquad
u_y = {{\sin\theta}\over{\sinh\mu}}.
\tag 3
$$
Equation~(1) becomes
$$
{{\partial\hfil}\over{\partial x}}
\left({{\cos\theta}\over{\cosh\mu}}\right)
+
{{\partial\hfil}\over{\partial y}}
\left({{\sin\theta}\over{\cosh\mu}}\right)
= 0,
$$
and the identity $(u_x)_y = (u_y)_x$ becomes
$$
{{\partial\hfil}\over{\partial x}}
\left({{\sin\theta}\over{\sinh\mu}}\right)
-
{{\partial\hfil}\over{\partial y}}
\left({{\cos\theta}\over{\sinh\mu}}\right)
= 0.
$$
Solving this latter pair of equations for the partials of~$\theta$ yields the
formulae
$$
\aligned
\sinh 2\mu\>\theta_x &=\phantom{-}\sin(2\theta)\,\mu_x
-\bigl(\cosh 2\mu+\cos(2\theta)\bigr)\,\mu_y,\\
\sinh 2\mu\>\theta_y &=-\sin(2\theta)\,\mu_y
+\bigl(\cosh 2\mu-\cos(2\theta)\bigr)\,\mu_x .\\
\endaligned
\tag 4
$$

Differentiating the equations~(4) and equating mixed partials of~$\theta$
yields the relation
$$
(\mu_{xx}-\mu_{yy})\,\cos(2\theta)
+2\mu_{xy}\,\sin(2\theta) = (\mu_{xx}+\mu_{yy})\,\cosh(2\mu)
\tag 5
$$
which evidently would determine $\theta$ modulo $\pi$ up to a choice of at most
two solutions unless the coefficients of $\cos(2\theta)$ and $\sin(2\theta)$
in~(5) were to vanish identically.  This observation motivates breaking the
analysis into two cases.

\subhead 1.2.  Degenerate cases \endsubhead
Suppose that~$\mu$ satisfies the condition that~(5) is an identity as an
equation for~$u$.  This is equivalent to requiring that $\mu$ satisfy
$$
\mu_{xx} = \mu_{xy} =\mu_{yy} = 0.
$$
In other words, $\mu$ is a linear function of $x$ and $y$:
$$
\mu(x,y) = ax+by +c
$$
for some constants $a$, $b$, and $c$ (not all zero).

\subsubhead 1.2.1.  Planes \endsubsubhead
If $a = b = 0$, then $\mu = c$ is a constant, so the equations~(4) imply that
$\theta$ is a constant and hence, by~(3), that $u$ is a linear function
of~$x$~and~$y$.  Of course, this is the first case of the theorem.

\subsubhead 1.2.2.  Scherk's Surfaces \endsubsubhead
On the other hand, if at least one of $a$ or $b$ is non-zero, then by rotating,
translating, and scaling if necessary, I can reduce to the case $\mu=x$.  Since
$\mu>0$, it follows that $\cD$ is a subdomain in the right half plane~$x>0$, so
I may as well let~$\cD$ be this entire half-plane.  The equations~(4) can now
be integrated, with the result that there is a constant~$\psi$ so that
$$
\theta(x,y) = -\tan^{-1}\bigl(\tanh(x)\,\tan(y+\psi)\bigr).
$$
Now integrating~(3), the corresponding 1-parameter family of inequivalent
solutions~$u$ is
$$
u(x,y) = \sinh^{-1}\left({{\cos (y+\psi)}\over{\sinh x}}\right).
$$
Thus, the graph of $u$ is one-half of one of the classical Scherk surfaces,
defined in $xyz$-space by the equation
$$
\sinh x\,\sinh z = \cos (y+\psi).
$$
Note that the corresponding area form is $dA = \coth x\, dx\w dy$ and is
invariant under translation in the $y$-direction.  This constitutes the second
case of the theorem.

\subhead 1.3.  Higher compatibility conditions \endsubhead
Now, I want to turn to the case where equation~(5) is not an identity.  In
order for there to be any solutions~$\theta$, the quantity
$$
P = \bigl(\mu_{xx}-\mu_{yy}\bigr)^2
+\bigl(2\mu_{xy}\bigr)^2
-\cosh^2(2\mu)\,\bigl(\mu_{xx}+\mu_{yy}\bigr)^2
\tag 6
$$
must be non-negative.  Set $\Delta = (\mu_{xx}-\mu_{yy})^2
+(2\mu_{xy})^2$.  If $\Delta=0$, then $P\ge0$ implies that
$\mu_{xx}=\mu_{xy}=\mu_{yy}$, which has already been studied.  Thus, from now
on I assume that $\Delta>0$.  The two solutions of~(5) are found to be
$$
\pmatrix \cos2\theta \\ \sin2\theta \endpmatrix
= \pmatrix \mu_{xx}-\mu_{yy}\\ 2\mu_{xy}\endpmatrix
{{(\mu_{xx}+\mu_{yy})\cosh 2\mu}\over{\Delta}}
\pm
\pmatrix-2\mu_{xy}\\ \mu_{xx}-\mu_{yy} \endpmatrix
{{\sqrt{P}}\over{\Delta}}
\tag 7
$$
{}From now on, I am going to assume that $P>0$, since otherwise, there is at
most one solution~$\theta$ modulo integral multiples of~$\pi$ and,
consequently, at most one minimal surface (up to equivalence) whose area
function is~$F$.  Substituting each of these two solutions into the
equations~(4) yields four equations of third order for the function~$\mu$.
These equations are fairly ugly and I won't write them out.  What is remarkable
is that making a change of dependent variable by setting $C = \cosh 2\mu$
results in equations for~$C$ which are quite simple:
$$
\aligned
C_{xxx} &=  (C_x\,C_{xx}-C_x\,C_{yy}+C_y\,C_{xy})/C\\
C_{xxy} &=  (C_x\,C_{xy})/C\\
C_{xyy} &=  (C_y\,C_{xy})/C\\
C_{yyy} &=  (C_y\,C_{yy}-C_y\,C_{xx}+C_x\,C_{xy})/C\,.\\
\endaligned
\tag 8
$$

Now, in general, one doesn't expect an overdetermined system such as~(8) to
have any solutions at all, but this system turns out to be involutive (as is
easily checked) and hence at any point $(x_0,y_0)$, one can freely specify the
values of $C (\not=0)$, $C_x$, $C_y$, $C_{xx}$, $C_{xy}$, and $C_{yy}$ and get
a unique solution to the above system of equations.  It follows that there is a
six parameter family of solutions.  Moreover, these equations are seen to have
the following first integrals:  For any solution of~(8), there exist constants
$a_1$, $a_2$, $a_3$ so that
$$
\aligned
C_{xy} &= a_1\,C\\
C_{xx}-C_{yy} &= a_2\,C\\
C\,(C_{xx}+C_{yy}) - C_x^2 - C_y^2 &= a_3\\
\endaligned
\tag 9
$$
Moreover, for any choice of the constants $a_i$, any solution of~(9) also
solves~(8).

To write out the general solution of~(9), I will take advantage of the
invariance of these equations under the group generated by translations,
rotations, and scaling in the $xy$-variables.  Rotating the independent
variables $x$ and $y$ by an angle~$\phi$ rotates the pair $(a_1,a_2)$ by an
angle of~$2\phi$.  Thus, I can assume that $a_1=0$ and $a_2\ge0$.

\subsubhead 1.3.1. The Heli-Catenoids \endsubsubhead
If $(a_1,a_2)=(0,0)$, then the solutions of~(9) are of the form
$$
C(x,y) = a(x^2+y^2) + bx + cy + d
$$
for some constants $a$, $b$, $c$, and $d$.  It is easy to see that if $a=0$,
then $P=0$, a case already discarded.  Thus, $a\not=0$ and so, by translation
in the $xy$-plane, I can assume $b=c=0$.  By direct calculation, the condition
$P>0$ is now seen to imply that $d^2<1$.  Since $C=\cosh 2\mu\ge 1$, it then
follows that $a$ must be positive.  By scaling in the $xy$-plane, I can
therefore assume that $a=2$ and take $d=\cos2\phi$.  The corresponding formula
for $F$ is then found to be
$$
F_\phi(x,y) = \sqrt{{x^2+y^2+\cos^2\phi}\over{x^2+y^2-\sin^2\phi}}
\tag 10
$$
where $0<\phi<\pi/2$ and of course, $\cD$ is defined by the inequality
$x^2+y^2>\sin^2\phi$.

A moment's thought leads to a guess as to what the corresponding minimal
surfaces are:  They are the surfaces in the isometric family of associated
minimal surfaces containing the helicoid ($\phi=0$) and the catenoid
($\phi=\pi/2$), i.e.,  the so-called ``heli-catenoids''.  The two inequivalent
solutions~$u$ are got by simply following the upper and lower halves of the
catenoid through the deformation in the associated family.  A direct
calculation, which I won't do here, verifies this guess.

\subsubhead 1.3.2. Doubly Periodic Solutions \endsubsubhead
If $a_1=0$ but $a_2>0$, then by scaling in the $xy$-plane, I may further assume
that $a_2=1$.  The equation $C_{xy}=0$ then implies that $C(x,y) =
\alpha(x)+\beta(y)$ for some functions $\alpha$ and $\beta$ of single
variables.  The equation $C_{xx}-C_{yy}=C$ becomes
$$
\alpha''(x)- \alpha(x) = \beta''(y)+\beta(y).
$$
It follows that each side of this equation must represent a constant function.
By subtracting this constant from $\beta$ and adding it to $\alpha$, I can
arrange that $\alpha''(x)=\alpha(x)$ and $\beta''(y) = -\beta(y)$.  The third
equation of~(9) is now an identity.

The conditions that $P>0$ and that $C\ge 1$ on~$\cD$ then imply that
appropriately translating in~$x$~and~$y$ will put the solution in the form
$$
C(x,y) = a\,\cosh x + c\,\cos y
$$
for some constants $a$ and $c$ satisfying the inequalities $|a{-}c|<1<a{+}c$.
(This describes a semi-infinite strip in the first quadrant of the $ac$-plane.)
 The corresponding domain~$\cD_{a,c}$, i.e.,  where $C >1$, is the $xy$-plane
minus a sequence of congruent `ovals' centered on the points
$(x,y)=\bigl(0,(2k+1)\pi\bigr)$. (The size of these `ovals' depends on the
values of $a$ and $c$.)

The function~$F$ then takes the form
$$
F_{a,c}(x,y) = \sqrt{{a\,\cosh x+c\,\cos y+1}\over{a\,\cosh x+c\,\cos y-1}}\>.
\tag 11
$$

To get the function~$\theta$, or more precisely, the functions~$\cos\theta$ and
$\sin\theta$ (which are all that are actually needed), one must now solve~(5),
which has now become
$$\displaylines{
\bigl(c^2+2ac\,\cosh x\,\cos y+a^2-1\bigr)\,\cos 2\theta
+2ac\,\sinh x \,\sin y \,\sin 2\theta\hfill\cr
\hfill=a\,(a^2-c^2-1)\,\cosh x +c\,(a^2-c^2+1)\,\cos y.\cr}
$$
The analysis of this equation is simplified by introducing the smooth
surface~$\Sigma_{a,c}$ in~$\bbR^3$ defined by the relation
$$
z^2 = a\,\cosh x + c\,\cos y - 1.
$$
Away from its intersection with the plane $z=0$, this surface is two smooth
copies of $\cD_{a,c}$, and may be pictured roughly as two copies of the
$xy$-plane joined by a countable number of ``throats'' evenly spaced along the
$y$-axis.  One computes that $\cos 2\theta$ and $\sin 2\theta$ can be smoothly
and globally defined as functions on~$\Sigma_{a,c}$ so as to satisfy the above
equation by means of the formulae
$$\pmatrix \cos 2\theta\\ \sin 2\theta\endpmatrix
= \frac{1}{A^2+B^2}
\pmatrix A & -B\\ B & \phantom{-}A\endpmatrix
\pmatrix E\\ Qz\endpmatrix
$$
where $A$, $B$, $E$ and $Q$ are smooth functions on the $xy$-plane defined by
the formulae
$$
\align
A &= c^2+2ac\,\cosh x\,\cos y+a^2-1\\
B &= 2ac\,\sinh x \,\sin y\\
E &= a\,(a^2-c^2-1)\,\cosh x +c\,(a^2-c^2+1)\,\cos y\\
Q &= \bigl((1-(a-c)^2)((a+c)^2-1)(a\,\cosh x + c\,\cos y + 1)\bigr)^{1\over2}.
\endalign
$$
Note that $A^2+B^2$ and~$Q$ are strictly positive on $\Sigma_{a,c}$, so these
formulae do, indeed define $\cos2\theta$ and $\sin2\theta$ as smooth functions
on~$\Sigma_{a,c}$.  By construction, these two functions are periodic in the
$y$-direction, of period~$2\pi$, and have the same values at the points
$(x,y,z)$ and $(-x,-y,z)$.

Now, I claim that it is possible to define $\cos\theta$ and $\sin\theta$ as
smooth functions on~$\Sigma_{a,c}$.  To see this, it suffices to show that, as
a multi-valued function on~$\Sigma_{a,c}$, the function~$\theta$ changes by an
integral multiple of~$2\pi$ as one traverses any closed loop on~$\Sigma_{a,c}$.
 Because of the $y$-periodicity, it suffices to check that this happens for two
specific loops, say $\sigma$ and $\gamma$, where $\sigma$ is any one of the
loops on the surface cut out by the plane $x=0$ and $\gamma$ is any one of the
loops on the surface cut out by the plane $z=0$.  For the first loop~$\sigma$,
it is not hard to see that $\theta$ does not change value at all when one goes
around it.  The second loop is somewhat more difficult to understand directly,
so I will replace it by a homotopic loop~$\gamma_R$.   Fix a large number
$R\gg0$ and consider the curve~$\gamma_R$ on the upper sheet of $\Sigma_{a,c}$
which lies over the boundary of the rectangle $|x|\le R$ and $0\le y \le 2\pi$.
 The corners $P_1$, $P_2$, $P_3$, and $P_4$ lie above the points $(-R,0)$,
$(R,0)$, $(R,2\pi)$, and $(-R,2\pi)$, respectively.  Now, dividing the original
equation above by $e^x$ and assuming $x=R$ (where $R$ is large), one sees that
$$
\cos\bigl(2\theta(R,y,z)-y\bigr) \approx \frac{a^2-c^2-1}{2c},
$$
while dividing the original equation by $e^{-x}$ and setting $x=-R$, one sees
that
$$
\cos\bigl(2\theta(-R,y,z)+y\bigr) \approx \frac{a^2-c^2-1}{2c},
$$
where the approximations are uniformly close in~$R$.  Since $|a^2-c^2-1|<|2c|$,
it follows that when $x\gg0$ the function~$\theta(x,y,z)$ is approximately
$\frac{1}{2} y+\theta_0$ where $\theta_0$ is a constant and when $x\ll0$ the
good approximation is $-\frac{1}{2} y + \theta_0$.  Using this and the evenness
of $\cos2\theta$ and $\sin2\theta$, one sees that, as one goes counterclockwise
around the path~$\gamma_R$, the multivalued function~$\theta(x,y,z)$ increases
by~$2\pi$.  This, plus the $2\pi$-periodicity in~$y$ of the functions $A$, $B$,
$E$, and~$Q$, shows that $\cos\theta$ and $\sin\theta$ can be smoothly defined
on the entire surface~$\Sigma_{a,c}$.  Moreover, they are periodic in~$y$ of
period~$4\pi$ and satisfy
$$
\align
\cos \theta(x,y+2\pi,z) &= -\cos\theta(x,y,z)\\
\sin \theta(x,y+2\pi,z) &= -\sin\theta(x,y,z)
\endalign
$$

Now, the desired function~$u$ on $\cD_{a,c}$ satisfies
$$
du = {{\cos\theta\,dx + \sin\theta\,dy}
\over{\sqrt{{\ts{1\over2}(a\,\cosh x + c\,\cos y-1)}}}},
\tag 12
$$
where the 1-form on the right hand side of~(12), say~$\zeta$, can be shown to
be globally and smoothly defined on~$\Sigma_{a,c}$.  Note that when $|x|$ is
large the 1-form $\zeta$ is quite small, so that the desired function~$u$ is
nearly constant for $|x|$ large, i.e., its graph is asymptotic to a plane.

Set
$$
\lambda_{a,c} = \int_\sigma \zeta.
$$
It is not hard to show that $\lambda_{a,c}$ is non-zero (which sign it has
depends on which choice of $(\cos\theta,\sin\theta)$ was made).  Let $L_{a,c}$
be the lattice in~$\bbR^3$ generated by $(0,4\pi,0)$ and $(0,0,\lambda_{a,c})$.
 Then the image of $\Sigma_{a,c}$ in $\bbR^3/L_{a,c}$ is easily seen to be a
torus~$T_{a,c}$ with four punctures.  The 1-forms $dx$, $dy$, and $\zeta$ are
well-defined on~$T_{a,c}$ and the period vectors of the vector-valued
1-form~$\xi$ on~$T_{a,c}$ defined by
$$
\xi = \pmatrix dx \\ dy\\ \zeta\endpmatrix
$$
lie in~$L_{a,c}$ by construction.  Thus there is a smooth
map~$X:T_{a,c}\to\bbR^3/L_{a,c}$ so that~$dX = \xi$.  By its very construction,
the image of~$X$ is a minimal surface in~$\bbR^3/L_{a,c}$ endowed with its
quotient flat structure as an abelian group.  Moreover the information derived
so far shows that this minimal surface is complete and of finite total
curvature.  In fact, endowing~$T_{a,c}$ with the conformal structure and
orientation induced on it by its immersion into~$\bbR^3/L_{a,c}$, one sees that
it is biholomorphic to a Riemann surface of genus~1 punctured at four points.

This analysis can be carried further to show that the Weierstrass
representation of this minimal surface is such that the completion of~$T_{a,c}$
to a curve of genus~1 can be identified with $\bbC/\Lambda$, where
$\Lambda\subset\bbC$ is a rank~2 lattice generated by a real and a purely
imaginary period.  It is not hard to see that $\zeta =
\text{Re}\bigl(2\,dw\bigr)$ where $dw$ is a holomorphic differential
on~$\bbC/\Lambda$.   Letting $g\colon \bbC/\Lambda\to S^2\simeq \bbC\bbP^1$ be
the meromorphic function which represents the Gauss map of the minimal surface,
one sees that
$$
\pmatrix dx \\ dy \\ \zeta\endpmatrix
=\text{Re}
\pmatrix
\hfil \bigl(g-g^{-1}\bigr)\,dw\\
i\,\bigl(g+g^{-1}\bigr)\,dw\\
2\,dw\endpmatrix .
$$
Note that $g$ has two simple poles at two of the half-lattice points and two
simple zeros at the other two half-lattice points. Thus, these minimal surfaces
are of the type studied by Karcher~\cite{Ka}.

To recover the two minimal graphs, let $S_{a,c}\subset\bbR^3$ denote the
inverse image of the minimal surface
$X\bigl(T_{a,c}\bigr)\subset\bbR^3/L_{a,c}$.  After a suitable vertical
translation, the part of~$S_{a,c}$ which lies between the two planes $z=0$ and
$z = \lambda_{a,c}$ will double cover~$\cD_{a,c}$ with a fold along the `oval'
boundaries.  These two halves will then be the two inequivalent minimal graphs
which induce the same area function.

This completes the proof of Theorem~1.

\head 2.  Special extremals for isosystolic metrics  \endhead

\subhead 2.1. The isosystolic ratio \endsubhead
I now want to describe a second problem of this sort which is motivated by
Professor Calabi's work on extremal isosystolic metrics on surfaces.  For any
Riemannian metric~$g$ on an $n$-manifold~$M$, Gromov~\cite{Gr} has defined the
{\it systole\/}~$\text{Sys}(M,g)$, to be the length of the shortest closed
curve~$\gamma$ in~$M$ which is not null-homotopic.  A manifold~$M$ is said to
be {\it essential\/} if there is a constant~$c(M)>0$ so that, for all
metrics~$g$ on~$M$,
$$
\text{Sys}(M,g) \le c(M)\, \bigl(\text{Vol}(M,g)\bigr)^{1/n}.
$$
The smallest possible value of~$c(M)$ is called the {\it minimum isosystolic
ratio} for~$M$.

While it is true that every compact surface other than the sphere is essential,
there are only three surfaces for which the extremal isosystolic ratio is known
exactly, the real projective plane, the torus, and the Klein bottle.  Professor
Calabi~\cite{Ca2} has set himself the task of computing or closely estimating
the isosystolic ratio for other surfaces and, along the way, has turned up
quite an interesting problem in the calculus of variations.

\subhead 2.2. Professor Calabi's work \endsubhead
Let me explain how this calculus of variations problem arises.  In searching
for an extremal isosystolic metric, Calabi has narrowed the search to a class
of (singular) generalized Riemannian metrics on a given surface~$M$ which have
the property that, aside from a small singular set, the surface can be covered
by ``tight bands'' of shortest non-trivial geodesics, i.e., geodesics which are
systole-long.  Outside of this singular set, the surface can be broken up into
what Calabi calls ``$k$-regular domains''~$U_k$ through which $k$ such bands
pass.

Calabi first shows that, in any 2-regular domain~$U_2$ for an extremal
isosystolic metric~$g$, the two bands must meet at right angles and that the
metric~$g$ must be flat there.  He does this by showing that, otherwise, he can
deform the metric in a neighborhood of such a point so as not to decrease the
lengths of the shortest geodesics while at the same time decreasing the area,
thus violating the assumption that $g$ was extremal.

Calabi then turns his attention to the problem of writing down the possible
metrics in a 3-regular domain~$U_3$.  He shows that each of the three bands
$B_i\cap U_3$ is calibrated by a closed 1-form of unit length $\omega^i$, i.e.,
a closed 1-form of length one with respect to $g$ which restricts to each of
the geodesics in the band~$B_i$ to be the differential of arc-length.  Of
course, these three 1-forms can be written locally in the form $\omega^i =
dx^i$.  We can take the first two as coordinates, say $x^1=x$ and $x^2=y$, and
the consider the third as a function~$z=F(x,y)$ of the first two.  The
condition that the metric $g$ make each of~$dx$, $dy$, and $dz$ be of unit
length is easily seen to be the condition that
$$
\bigl|\alpha\,dx + \beta\,dy\bigr|^2_g = \alpha^2 + \beta^2 +
2f(x,y)\,\alpha\beta
$$
where $f$ is given by the expression
$$
f(x,y) = {{1-F_x^2-F_y^2}\over{2F_xF_y}}
$$
and the corresponding area form is
$$
dA_g = {{\bigl|2\,z_x\,z_y\bigr|}
\over{\sqrt{(2\,z_x\,z_y)^2-(z_x^2+z_y^2-1)^2}}}\,dx\w dy.
$$
It follows that, if the metric $g$ is to be extremal isosystolic, then the
function~$z=F(x,y)$ must be a solution of the Euler-Lagrange equation of the
functional
$$
\cF(z) = \int_\cD {{2\,z_x\,z_y}
\over{\sqrt{(2\,z_x\,z_y)^2-(z_x^2+z_y^2-1)^2}}}\,dx\w dy\,,
$$
at least in the region where the three bands are tight.

It is not difficult to verify that the Euler-Lagrange equation for this
functional is elliptic at solutions which satisfy the inequalities
$z_x^2{+}z_y^2>1$ and $|z_x^2{-}z_y^2|<1$.  Now, because the Lagrangian
$$
L(x,y,z,p,q) = {{2\,p\,q}
\over{\sqrt{(2\,p\,q)^2-(p^2+q^2-1)^2}}}
$$
does not involve $x$, $y$, or $z$ explicitly, the Legendre transform of its
Euler-Lagrange equation becomes an equation for a function of~$p$ and~$q$
which, on the region~$\cE$ in $pq$-space defined by the inequalities
$p^2{+}q^2>1$ and $|p^2{-}q^2|<1$, is a second-order, linear, elliptic {\smc
pde} of the form $\Delta \Phi = 0$ where $\Delta$ is the Laplacian of an
appropriate metric on~$\cE$.  Of course, this equation can be reduced to
Laplace's equation once one finds a complex coordinate~$\zeta:\cE\to\bbC$ which
uniformizes the induced conformal structure.  (This is the analog of
introducing conformal coordinates on a minimal surface in space and deriving
the Weierstrass representation for minimal surfaces:  The equation is reduced
to a linear one by a change of coordinates.)  However, such a coordinate has
proved rather difficult to find.  In fact, no solutions of the corresponding
Euler-Lagrange equation other than those which lead to a flat metric $g$ are
known, and these trivial solutions are useless for finding a uniformizing
coordinate on~$\cE$.

\subhead 2.3. The prescribed Lagrangian problem \endsubhead
To make further progress, at least one non-trivial solution to the
Euler-Lagrange equation needs to be found.  The usual approach for finding
non-trivial solutions would be to look for special solutions, particularly
solutions which are invariant under some continuous group of symmetries of the
equation.  Unfortunately, since the Lagrangian has relatively few symmetries,
one cannot expect solutions with continuous symmetries and, in fact, no
non-trivial solutions invariant under a continuous symmetry of the equation are
known.

However, returning to the original problem, Calabi realized that understanding
the extremal isosystolic metrics on 4-regular domains involves finding two
solutions $z$ and $w$ to the Euler-Lagrange equation which induce the same
metric~$g$, i.e., which satisfy
$$
{{2\,z_x\,z_y}
\over{\sqrt{(2\,z_x\,z_y)^2-(z_x^2+z_y^2-1)^2}}}
=\pm \
{{2\,w_x\,w_y}
\over{\sqrt{(2\,w_x\,w_y)^2-(w_x^2+w_y^2-1)^2}}}\,.
$$
More generally, in a $k$-regular region (where $k\ge3$), one expects that for
the area form $dA_g = F(x,y)\,dx \w dy$ there should be $k{-}2$ distinct
solutions~$z^i$ to the Euler-Lagrange equation which also satisfy
$$
{{2\,z^i_x\,z^i_y}
\over{\sqrt{(2\,z^i_x\,z^i_y)^2-\bigl((z^i_x)^2+(z^i_y)^2-1\bigr)^2}}}
 = \pm\,F.
$$
Thus, this is a prescribed density problem of the sort under discussion.
Again, since the Lagrangian does not contain~$z$ explicitly and because it is
an even function of $(p,q)$, it follows that, for any solution $z$ to the
problem, every function of the form~$\pm z+c$ is also a solution.  However,
such solutions are to be regarded as equivalent because they give rise to the
same local foliation by geodesics.  I can now state the main result:

\proclaim{Theorem~2}For any non-constant function $F$ on a domain $\cD$ in the
plane, there are at most two inequivalent solutions to the Euler-Lagrange
equation for~$L$ which satisfy $L=\pm F$.  Moreover, up to translation and
dilation, the space of those~$F$ for which there are two inequivalent solutions
to this Euler-Lagrange satisfying $L=\pm F$ is of dimension at most~3.
\endproclaim

This theorem has the following corollary:

\proclaim{Corollary} For any extremal isosystolic metric~$g$, the Gauss
curvature vanishes in all $k$-regular domains with $k\ge 5$.  Moreover, there
is a family~$\cF$ of metrics, depending on at most~3 parameters, so that in a
$4$-regular domain an extremal isosystolic metric~$g$ is locally isometric to
one of the metrics in this family.
\endproclaim

I conjecture that, in fact, the space of~$F$ for which there are two
inequivalent solutions consists only of the constants.  Of course, if this
conjecture is true, then the corollary about local flatness can be sharpened to
include the 4-regular domains as well.  This would leave only the 3-regular
domains as unknown open subsets of the extremal isosystolic metrics.

The proof of Theorem~2 rests on an analysis of the compatibility conditions for
the unknown function $F$ with $L$ given as above.  These calculations are quite
long and were carried out on a computer, using the symbolic manipulation
program {\smc maple}.  For this reason, at a certain point, I am only going to
indicate how the calculations are to be carried out.

\subhead 2.4. Reduction to an equation\endsubhead
Here is the setting of the problem in the elliptic case:  For a given function
$F\ge1$ on a domain $\cD$ in the $xy$-plane, find all of the inequivalent
solutions~$z$ on~$\cD$ which satisfy the Euler-Lagrange equations for the
functional
$$
\cL(z) = \int_\cD {{2\,z_x\,z_y}
\over{\sqrt{(2\,z_x\,z_y)^2-(z_x^2+z_y^2-1)^2}}}\,dx\w dy
$$
and which also satisfy the first order equation
$$
{{2\,z_x\,z_y}
\over{\sqrt{(2\,z_x\,z_y)^2-(z_x^2+z_y^2-1)^2}}}
=F.
$$

Setting $F = 1/\sin(2\phi)$ where $\phi$ is a function on~$\cD$ with values in
the range $0<\phi < \pi/4$, the above first order equation can be put in the
form
$$
(2\,\cos2\phi\,z_x\,z_y)^2 = (z_x^2+z_y^2-1)^2.
$$
Considering the elliptic range and reflecting in $x$ and/or $y$ if necessary, I
can assume that $z_x$ and $z_y$ are positive and satisfy $z_x^2{+}z_y^2>1$ and
$|z_x^2{-}z_y^2|<1$.  Then taking the positive square root of both sides yields
$$
z_x^2 - 2\,\cos2\phi\,z_x\,z_y + z_y^2=1.
$$
Thus, $(z_x,z_y)$ is determined to lie on an ellipse which passes through the
points $(1,0)$ and $(0,1)$.  From this and the ellipticity assumption, it is
easy to see that there must be a function~$\theta$ satisfying
$|\theta|<\frac{1}{2}\pi-\phi$ so that
$$
z_x = {{\cos (\theta-\phi)}\over{\sin 2\phi}}\,,
\qquad\text{and}\qquad
z_y = {{\cos (\theta+\phi)}\over{\sin 2\phi}}\,.
\tag 1
$$
Now, the Euler-Lagrange equation is equivalent to the condition that the
\hbox{1-form} $\psi$ be closed where
$$
\psi =
{{(z_x^4-z_y^4-2z_x^2+1)z_x\,dx-(z_y^4-z_x^4-2z_y^2+1)z_y\,dy}
\over{\left((2\,z_x\,z_y)^2-(z_x^2+z_y^2-1)^2\right)^{3\over2}}}.
$$
Substituting the formulae for $z_x$ and $z_y$ into the equations
$d\psi = d(z_x\,dx+z_y\,dy)=0$ and expanding allows one to solve for the
partials of the unknown function $\theta$ in terms of the known function~$\phi$
and its partials.  The actual formulae are rather complicated, but for what I
am going to describe, it suffices to know that the result is an equation of the
form
$$
d(2\theta) = \cos 2\theta\,\omega_1+\sin 2\theta\,\omega_2+\omega_3
\tag 2
$$
where $\omega_1$, $\omega_2$, and $\omega_3$ are certain 1-forms which are
linear combinations of $dx$ and $dy$ with coefficients given by some (rather
ugly) rational expressions in $\sin \phi$, $\cos \phi$, $\phi_x$, and $\phi_y$
(and linear in the latter two variables).

Differentiating~(2) yields
$$
0 = d\bigl(d(2\theta)\bigr)
= (A_1\,\cos 2\theta +A_2\,\sin 2\theta+A_3)\, dx\w dy
$$
where $A_1$, $A_2$, and $A_3$ are some (also rather ugly) rational expressions
in $\sin \phi$, $\cos \phi$, $\phi_x$, $\phi_y$, $\phi_{xx}$, $\phi_{xy}$, and
$\phi_{yy}$ (and linear in the last three variables).  It follows that any
solution $\theta$ must satisfy the relation
$$
A_1\,\cos 2\theta +A_2\,\sin 2\theta+A_3 = 0.
\tag 3
$$
(This is the analog of equation~(5) in \S1.1.)

Just as in the case of minimal surfaces, there are now two cases to consider:
The first is that $\phi$ might satisfy the second order quasi-linear {\smc pde}
system
$$
A_1 \equiv A_2 \equiv A_3 \equiv 0.
$$
Just as in the case of minimal surfaces, these three equations can be solved
for $\phi_{xx}$, $\phi_{xy}$, and $\phi_{yy}$, so that this system can be
written in the form
$$
\phi_{xx} = F_1(\phi,\phi_x,\phi_y),\quad
\phi_{xy} = F_2(\phi,\phi_x,\phi_y),\quad
\phi_{yy} = F_3(\phi,\phi_x,\phi_y).
$$
However, unlike the analogous system for minimal surfaces treated in \S1.2,
this system is not involutive.  In fact, an analysis of these equations using
{\smc maple} shows that the only solutions of this system are given by
setting~$\phi$ equal to a constant.  Of course, in this case, $\theta$ must
also be constant and hence $z$ must be a linear function of $x$ and $y$,
yielding the known case where there are inequivalent solutions.

The second case is that not all of the $A_i$ vanish, so that (3) becomes a
non-trivial equation for~$\theta$.  In fact, clearly there can be at most two
solutions~$\theta$ of~(3) which lie in the elliptic range
$|\theta|<\frac{1}{2}\pi-\phi$.  The condition that there actually be two
distinct solutions is a set of second order differential inequalities on~$\phi$
that it is not helpful to write out here.

Supposing that there are two distinct solutions, say $\theta^\pm$, they can be
written out explicitly in terms of $\phi$ and its first and second derivatives
in a formula analogous to equation~(7) of~\S1.3.  These expressions are rather
large and I won't write them out here.  It remains to check whether these two
solutions of~(3) actually satisfy~(2), i.e.,  that the following equations
hold:
$$
d(2\theta^\pm)
= \cos 2\theta^\pm\,\omega_1+\sin 2\theta^\pm\,\omega_2+\omega_3.
\tag 4
$$
Substituting the explicit formulae for $\theta^\pm$, the equations~(4) expand
to a system of four third-order equations for~$\phi$.  In fact, these equations
can be solved in terms of the third derivatives in the form
$$
\aligned
\phi_{xxx} &= G_1(\phi,\phi_x,\phi_y,\phi_{xx},\phi_{xy},\phi_{yy}),\\
\phi_{xxy} &= G_2(\phi,\phi_x,\phi_y,\phi_{xx},\phi_{xy},\phi_{yy}),\\
\phi_{xyy} &= G_3(\phi,\phi_x,\phi_y,\phi_{xx},\phi_{xy},\phi_{yy}),\\
\phi_{yyy} &= G_4(\phi,\phi_x,\phi_y,\phi_{xx},\phi_{xy},\phi_{yy}).\\
\endaligned
\tag 5
$$
Unfortunately, the functions $G_i$ are gross and I have not been able to check
their integrability conditions.  It could well be that the only solutions are
got by setting~$\phi$ equal to a constant.  In any case, only for a solution
$\phi$ of this system (satisfying the appropriate second order differential
inequalities) do there exist two inequivalent solutions $z$ to the original
problem with prescribed area form $F\,dx\w dy = \csc 2\phi\,dx\w dy$.

What can be said with confidence is that there is at most a six parameter
family of solutions~$\phi$ of~(5) (which would be the case if and only if this
system were in involution, as was the analogous system~(8) of \S1.3).  These
solutions reduce to at most three after taking into account the fact that the
whole problem is invariant under translation in $x$ and $y$ and simultaneous
dilation in $x$, $y$, and $z$.  (This implies that for any solution $\phi(x,y)$
of (5), all of the functions $\phi(cx+a,cy+b)$ are also solutions of (5).)
This observation completes the proof of Theorem~2.

Of course, it is possible that a clever change of variables will simplify the
equations~(5) to something manageable, but I haven't been able to find one.  My
suspicion is that the only solutions are the constants, but it looks as though
this will be hard to prove.  In theory, checking the integrability conditions
of (5) is routine, but, using {\smc maple}, it takes more memory than any
machine to which I have access.

\head 3. Harmonic maps of constant energy density between spheres \endhead

\subhead 3.1. Minimal surfaces and harmonic maps \endsubhead
In this last part, I want to consider another case where the notion of
solutions of Euler-Lagrange equations with prescribed Lagrangian makes an
interesting appearance.  Many years ago, Professor Calabi~\cite{Ca1}classified
the minimal 2-spheres of constant Gaussian curvature in the standard Euclidean
spheres, showing that they were simply the already known Boruvka spheres.  This
work was carried further by Wallach~\cite{W}, who showed that any minimal
submanifold $M^k\subset S^n$ of constant positive sectional curvature is a
(piece of a) global isometric immersion $f\colon (S^k,r)\to S^n$ given by
harmonic polynomials of the appropriate degree.

Now, the condition of being a minimal submanifold of constant positive
sectional curvature is very strong.  Recall that when one endows an
$n$-dimensional minimal submanifold with the induced metric, the inclusion
mapping becomes a harmonic map which is also an isometric embedding.  In
particular, its energy density is a constant~$n$ times the volume form.  Thus,
minimal submanifolds furnish examples of harmonic maps with constant energy
density.  Thus, it makes sense to generalize the study of minimal isometric
embeddings to that of harmonic maps of constant energy density.  Both of these
subjects have received attention in the literature from a number of authors
\cite{Br}, \cite{DW}, \cite{DZ}, \cite{Ke}, \cite{L}, \cite{PU}, \cite{T}.
Note that, since harmonic maps are just solutions to the Euler-Lagrange
equations obtained by taking the energy density as the Lagrangian, this sort of
problem fits exactly into the sort being discussed.

In~\cite{Br} it is proved that if~$\Sigma$ is a surface endowed with a  metric
of constant Gauss curvature~$K$ and $f\colon\Sigma^2\to S^n$ is a harmonic map
whose energy density $|df|^2$ is constant, then $K\ge 0$.  Moreover, if $K>0$
then $f(\Sigma)$ is a piece of a Boruvka sphere.%
\footnote{This rules out there being any minimal surfaces in any~$S^n$ of
constant negative curvature, which is reminiscent of a very early result of
Professor Calabi's to the effect that the Poincar\'e disk cannot be
isometrically embedded as a complex curve (even locally) in any~$\bbC\bbP^n$.}
Recently, I have extended this result to all dimensions.

\proclaim{Theorem~3}If $(\Sigma^n,K)$ is a \rom(piece of a\rom) space form of
constant sectional curvature~$K>0$ and $f\colon\Sigma\to S^N$ is a harmonic
mapping with constant energy density, then $f$ extends to a global harmonic map
from~$S^n$ to~$S^N$ which is given by harmonic polynomials on~$S^n$ of the
appropriate degree.
\endproclaim

The proof, which is based on an interesting identity, will be outlined in the
following sections.  Because the notation used below has certain failings when
$n=2$, I am going to assume from now on that $n$ is at least~3.  In any case,
the theorem is already proved in~\cite{Br} in the case~$n=2$.

\subhead 3.2. A Notation for Representations of the Orthogonal Group\endsubhead
The proof of Theorem~3 will involve working with symmetric tensors of rather
high order on a space form~$\Sigma^n$ of dimension~$n$.  For this reason, it is
important to establish an effective notation for the irreducible components of
the representations of~$\SO(n)$ on the symmetric powers of~$\bbR^n$ and various
bilinear pairings between them.

Let $x^1,\ldots, x^n$ denote the standard linear coordinate functions
on~$\bbR^n$.  Let $\cP$ denote the ring of polynomials in the variables~$x^i$
and let $\cH_d\subset \cP$ denote the vector space of harmonic polynomials of
degree~$d$.  It is known that $\cH_d$ is an irreducible representation of
$\SO(n)$.  Moreover, the symmetric tensor power~$S^d\bigl(\cH_1\bigr)$
of~$\cH_1\simeq\bbR^n$, regarded as the space of all homogeneous polynomials of
degree~$d$ on~$\bbR^n$, obeys the inductive decomposition rule
$$
S^d\bigl(\cH_1\bigr) = \cH_d\oplus\ R\cdot S^{d-2}\bigl(\cH_1\bigr),
$$
where $R=(x^1)^2 +\cdots+(x^n)^2$, thus implying the complete irreducible
decomposition of this symmetric power as an $\SO(n)$-module.

For any $f\in\cH_d$ and any $\xi\in\cH_1$, the product $\xi f$ is a homogeneous
polynomial of degree~$d{+}1$ and hence can be decomposed into irreducible
components by the above rule.  In fact, it turns out to be best for
normalization purposes to define
$$
(n+2d-2)\>\xi f = f\vee\xi + R\,(f\cdot \xi)
$$
where $f\vee\xi$ lies in $\cH_{d+1}$ and $f\cdot\xi$ turns out to lie in
$\cH_{d-1}$.  In fact, $f\cdot\xi$ turns out to simply be the directional
derivative of $f$ in the direction~$\xi$ (regarded, via the metric, as a vector
in~$\bbR^n$).  This unambiguously defines two pairings,
${\vee}\colon\cH_d\times\cH_1\to\cH_{d+1}$ and
${\cdot}\colon\cH_d\times\cH_1\to\cH_{d-1}$, which will be used extensively.
For example, when $n=3$, one obtains
$$
x^1\vee x^1 = 2(x^1)^2 - (x^2)^2 - (x^3)^2
\qquad\text{and}\qquad
x^1\cdot x^1 = 1.
$$

The Lie algebra~$\euso(n)$ of $\SO(n)$ is isomorphic to
$\Lambda^2(\bbR^n)\simeq\Lambda^2(\cH_1)$, and the formula for the  derived
action of this algebra on the representation~$\cH_{d}$ (regardless of the value
of~$d$) is found to be
$$
(\alpha\w\beta)\ldot f = \alpha\,(f\cdot \beta) - \beta\,(f\cdot \alpha).
$$
Moreover, there is an $\SO(n)$-invariant quadratic form on~$\cH_d$ defined by
the rule
$$
\langle f,g\rangle = \frac{(-\Delta)^d(fg)}{ 2^d d! }
$$
where $\Delta$ is the (geometer's) Laplacian, i.e., the {\it negative} of the
trace of the Hessian.  Finally, one can define an $\SO(n)$-equivariant pairing
$\{,\}:\cH_d\times\cH_{d+1}\to\cH_1$ by the rule
$$
\{f,\,g\}\cdot\alpha = \langle f,\,g\cdot \alpha\rangle\,.
$$

A rather large number of identities attend these pairings and the inner
product.  For example, for all~$f\in\cH_d$ and $\alpha,\beta\in\cH_1$,
$$
\alignedat2
\bigl(f\vee\alpha\bigr)\cdot\beta - \bigl(f\vee\beta\bigr)\cdot\alpha
&=&(n{+}2d)\>(\alpha\w\beta)\ldot f,& \\
\bigl(f\cdot\alpha\bigr)\vee\beta - \bigl(f\cdot\beta\bigr)\vee\alpha
&=&\ -(n{+}2d-4)\>(\alpha\w\beta)\ldot f,& \\
\bigl(f\vee\alpha\bigr)\cdot\alpha - \bigl(f\cdot\alpha\bigr)\vee\alpha
&=&(n{+}2d-2)\>||\alpha||^2\,f.&
\endalignedat
$$
Moreover, for all $f\in\cH_{d}$ and $g\in\cH_{d+1}$, the formula
$$
\langle f\vee \alpha,\,g\rangle = (n+2d-2)\,\langle f,\,g\cdot\alpha\rangle
$$
shows that, up to certain well-defined scale factors, the
operators~$\vee\alpha$ and~$\cdot\alpha$ are adjoint with respect to these
inner products.  Listing all of the identities to be used below would result in
a rather large, incomprehensible table.  Since the proofs of these and other
like identities are straightforward, I will merely point out when they are used
and leave the verifications to the reader.

\subhead 3.3. Structure equations on a space form\endsubhead
Now consider a local orthonormal coframing on~$\bigl(\Sigma^n,ds^2\bigr)$, say
$$
ds^2 = (\omega_1)^2 + \cdots + (\omega_n)^2
$$
By the first structure equation of \'E.~Cartan, there exist unique
1-forms~$\varphi_{ij}=-\varphi_{ji}$ representing the Levi-Civita connection on
the domain of the coframing so that
$$
d\omega_i = -\varphi_{ij}\w\omega_j\,.
$$
(Note the use, here and below, of the summation convention, even though all of
the indices have been lowered.)  Since the metric $ds^2$ has constant sectional
curvature~$K$, the second structure structure equation of \'E.~Cartan becomes
$$
d\varphi_{ij} = - \varphi_{ik}\w\varphi_{kj} + K\,\omega_i\w\omega_j\,.
$$

Now, these equations can be written in vector form as follows.  Let $\omega$
and $\varphi$ denote the 1-forms with values in $\cH_1$ and
$\Lambda^2(\cH_1)\simeq \euso(n)$, respectively, defined by the formulae
$$
\omega = \omega_i\,\,x^i
\qquad\text{and}\qquad
\varphi = {\ts\frac{1}{2}}\,\varphi_{ij}\,x^i\w x^j.
$$
Then the first and second structure equations are expressible in the form
$$
d\omega = {} - \varphi\,\ldot\,\omega
\qquad\text{and}\qquad
d\varphi = {} - {\ts\frac{1}{2}}\bigl[\varphi,\varphi\bigr]
                  + {\ts\frac{K}{2}}\,\omega\w\omega.
$$

\subhead 3.4. An identity for eigenfunctions of the Laplacian \endsubhead
Now let $f$ be a smooth function on the domain of the coframing.  There exists
a unique $\cH_1$-valued function~$f^{(1)}$ on this domain so that
$$
df = f^{(1)}\cdot\omega.
$$
(The components of~$f^{(1)}$ are basically the components of~$\nabla f$ in the
orthonormal coframing.)  Now, it is not difficult to see that there is a
formula of the form
$$
df^{(1)} = {} -\varphi\,\ldot\,f^{(1)} - {\ts\frac{1}{n(n-2)}}\,\Delta
f\vee\omega
+ f^{(2)}\cdot \omega,
$$
where $\Delta f$ is the (geometer's) Laplacian and $f^{(2)}$ is a
$\cH_2$-valued function which, up to scale, gives the components of the
traceless part of the Hessian of~$f$ in this coframing.  This formula could be
supplemented by formulae for~$df^{(2)}$, etc, but these rapidly become
complicated for a general function.  However, for eigenfunctions of the Laplace
operator, there is a very nice sequence of formulae:

\proclaim{Lemma}  Let $f$ satisfy $\Delta f = \lambda f$ for some
constant~$\lambda$.  Then there exists a sequence of functions $f = f^{(0)}$,
$f^{(1)}$, $f^{(2)}$, $\ldots$, with~$f^{(m)}$ taking values in~$\cH_m$, so
that $df^{(0)} = f^{(1)}\cdot \omega$ and
$$
df^{(m)} = -\varphi\,\ldot\,f^{(m)} - b_{m-1}(\lambda)\,f^{(m-1)}\vee\omega
+ f^{(m+1)}\cdot \omega
\qquad\qquad (m\ge 1),
\tag 1
$$
where the sequence $b_m(\lambda)$ for $m\ge 0$ is given by the formula
$$
b_m(\lambda) = \frac{\lambda-m(n+m-1)K}{(n+2m)(n+2m-2)}.
\tag 2
$$
Moreover, the function~$f$ satisfies $f^{(m+1)}\equiv0$ if and only if it is
the restriction to~$\Sigma$ of a homogeneous harmonic polynomial of degree~$m$
\rom(when $\Sigma$ \rom(or its simply connected cover\rom) is given its usual
interpretation as a hypersphere in a vector space of dimension~$n{+}1$\rom).
\endproclaim

The proof of the formula part of the Lemma is a straightforward induction,
using the pairing identities listed above and a few others very like them.  I
will not reproduce it here.  Of course, the function~$f^{(m)}$ is, up to a
scale factor, the traceless part of the fully symmetrized $m$-th covariant
derivative of~$f$.

The formula~(1) can be partly interpreted as a formula relating certain first
order differential operators
$$
\nabla^\pm \colon C^\infty\bigl(S^d_0(T^*)\bigr)\longrightarrow
C^\infty\bigl(S^{d\pm 1}_0(T^*)\bigr).
$$
In this language, suitably interpreted, equation~(1) implies that if $f\in
C^\infty(\Sigma)$ satisfies $\nabla^{-}\nabla^+ f = \lambda f$, then
$$
\nabla^-\bigl(\nabla^+)^m f = b_{m-1}(\lambda)\, (\nabla^+)^{m-1} f.
$$
However, the Lemma contains more information than this.  Moreover, the
definitions of the operators~$\nabla^\pm$ depend on choices of various scale
factors which are notoriously inconsistent in the reference literature.  Thus,
I have chosen to avoid this way of explaining the proof.

Seeing the second part of the Lemma, about the identification of the
eigenfunctions which satisfy $f^{(m+1)}\equiv 0$, require a few remarks.
First, of all, notice that, by~(1), if $m$ is the lowest value for which
$f^{(m+1)}$ vanishes, then $b_m(\lambda) f^{(m)}$ must vanish identically, so,
by the minimality of $m$, it follows that $b_m(\lambda)$ must vanish, i.e.,
that
$$
\lambda = m(m+m-1)K,
$$
which is exactly the eigenvalue of the restriction to~$\Sigma$ of a harmonic
polynomial of degree~$m$ on the ambient vector space in which~$\Sigma$ is
embedded in the standard embedding.  Moreover, setting $f^{(m+1)}=0$ and
considering the remaining differential formulae for $f^{(0)},\ldots,f^{(m)}$,
one sees that this becomes a set of linear total differential equations which
is closed under exterior derivative.  It follows that the vector space of
solutions to this system is of dimension
$$
\dim \cH_0 + \dim\cH_1 +\cdots + \dim \cH_m
= \frac{(n+2m-1)(n+m-2)!}{m! (n-1)!}
$$
which, again, is exactly the dimension of the space of harmonic polynomials on
$\bbR^{n+1}$ of degree~$m$.  Of course, this information is not conclusive.
The remainder of the argument comes from the following formula, which is easily
derived from the differential identity of the lemma and various pairing
identities:
$$
\displaylines{
d\bigl(\{f^{(m)},f^{(m+1)}\}\cdot\> {*}\omega\bigr)\hfill\cr
\hfill= \left[(m+1)\langle f^{(m+1)},f^{(m+1)}\rangle
- b_m(\lambda) (n+m-2)(n+2m) \langle f^{(m)},f^{(m)}\rangle\right]\>{*}1.\cr
}
$$
(Here, ${*}$ is the Hodge star operator.)  Because of the $\SO(n)$-invariance
of all of the quantities involved, the $n$-forms on either side of this
equation are well-defined on~$\Sigma$, independent of any choice of local
coframe.  Thus, if $\Sigma$ is compact, then integrating both sides of this
equation over~$\Sigma$ yields
$$
\int_\Sigma \langle f^{(m+1)},f^{(m+1)}\rangle\>{*}1
=
\frac{b_m(\lambda) (m+n-2)(n+2m)}{(m+1)}
\int_\Sigma\langle f^{(m)},f^{(m)}\rangle\>{*}1.
$$
{}From this, it follows immediately that if $\lambda = m(m+n-1)K$, then any
global eigenfunction~$f$ of the Laplacian with this eigenvalue (in particular,
the restrictions of the harmonic polynomials of degree~$m$ on $\bbR^{n+1}$)
must satisfy $\langle f^{(m+1)},f^{(m+1)}\rangle\equiv0$ and hence
$f^{(m+1)}\equiv 0$.  The remainder of the Lemma now follows by dimension
count.

It is amusing to note that the above integral formula shows directly that there
cannot be any non-zero eigenfunctions of the Laplacian on a compact~$\Sigma$
for any value of $\lambda$ not of the form $m(m{+}n{-}1)K$.  The reason is
that, if $\lambda$ is not of this form, then the fact that $K>0$ implies that
there will be a first value of $m$ for which~$b_m(\lambda)$ is negative.  Of
course, this will imply that the integral on the left hand side of this
equation is not positive, forcing $f^{(m+1)}\equiv0$, which, as has been shown,
cannot happen for any non-zero solution of $\Delta f = \lambda f$ unless
$\lambda=m(m+n-1)K$.

\subhead 3.5. Harmonic maps of constant energy density \endsubhead
Now suppose that $f\colon \Sigma\to S^N$ is a harmonic mapping with constant
energy density.  I will use $(,)$ to denote the inner product on~$\bbR^{n+1}$
which renders $S^N$ the unit sphere.. Thus, $f$ can be regarded as an
$\bbR^{N+1}$-valued function which satisfies $(f,f) = 1$.

By the usual computation for harmonic maps to spheres, the condition that $f$
be harmonic as a map from $\Sigma$ to $S^N$ is
$$
\Delta f = ||df||^2 f
$$
where $||df||^2$ is the energy-density of the mapping and $\Delta$ is the
(geometer's) Laplacian.  Thus, the condition that the energy-density be
constant is equivalent to the condition that the components of~$f$ be
eigenfunctions of the Laplacian, all with the same constant eigenvalue, say,
$\lambda$.

Looked at another way, one can say that constructing a (local) harmonic map of
constant energy density from $\Sigma$ to some $S^N$ is equivalent to choosing
$N{+}1$ (local) solutions of the equation $\Delta f_a = \lambda f_a$ with the
property that the sum of their squares is a constant function.  Essentially,
this can be thought of as $N{+}1$ second order {\smc pde} for $N$ unknown
functions, i.e., a (mildly) overdetermined partial differential equation, the
sort that the theory of exterior differential systems~\cite{BCG} was designed
to treat.

However, in this case, the problem is more easily treated directly.  In the
formulae below, the notation $\langle\!(,)\!\rangle$ will denote the tensor
product of the inner products $\langle,\rangle$ on $\cH_d$ and $(,)$ on
$\bbR^{N+1}$.  Similarly, given a function $g$ with values in
$\bbR^{N+1}\otimes\cH_d$ and $h$ with values in $\bbR^{N+1}\otimes\cH_{d+1}$,
the expression $\{\!(g,h)\!\}$ will denote the function with values in $\cH_1$
which is got by applying the inner product $(,)$ and the pairing $\{,\}$
defined above in the obvious way.

Now, because $f$ satisfies the vector equation $\Delta f = \lambda f$, the
formulae derived in the last section apply directly with little modification to
yield, for all~$m\ge 0$,
$$
\displaylines{
d\bigl(\{\!(f^{(m)},f^{(m+1)})\!\}\cdot\> {*}\omega\bigr)\hfill\cr
\hfill= \left[(m+1)\langle\!(f^{(m+1)},f^{(m+1)})\!\rangle
- b_m(\lambda) (n+m-2)(n+2m)
\langle\!(f^{(m)},f^{(m)})\!\rangle\right]\>{*}1.\cr
}
$$
Now consider the identity $(f,f) = \langle\!(f^{(0)},f^{(0)})\!\rangle
\equiv1$.  By differentiation, this implies that
$\{\!(f^{(0)},f^{(1)})\!\}\equiv0$.  Substituted into the above formula, this
recovers the already established identity
$\langle\!(f^{(1)},f^{(1)})\!\rangle\equiv\lambda$.  However, now
differentiating this identity yields $\{\!(f^{(1)},f^{(2)})\!\} \equiv 0$.
Setting $m=1$ in the above relation now gives
$$
\langle\!(f^{(2)},f^{(2)})\!\rangle
=\frac{b_1(\lambda)\,n(n+4)}{3}\langle\!(f^{(1)},f^{(1)})\!\rangle
=\frac{\lambda b_1(\lambda)\,n(n+4)}{3},
$$
so that, in particular, $\langle\!(f^{(2)},f^{(2)})\!\rangle$ is also constant.

Continuing now by induction, one sees that starting with the knowledge that
$\langle\!(f^{(m)},f^{(m)})\!\rangle$ is a constant function, say,
$\langle\!(f^{(m)},f^{(m)})\!\rangle = A_m(\lambda)\ge 0$, one can
differentiate, yielding $\{\!(f^{(m)},f^{(m+1)})\!\}\equiv 0$.  Substituting
this into the formula above yields
$$
\align
\langle\!(f^{(m+1)},f^{(m+1)})\!\rangle
&=\frac{b_m(\lambda) (m+n-2)(n+2m)}{(m+1)}
\langle\!(f^{(m)},f^{(m)})\!\rangle\\
&=\frac{b_m(\lambda) (m+n-2)(n+2m)}{(m+1)} A_m(\lambda).
\endalign
$$
Thus, each of the functions~$\langle\!(f^{(m)},f^{(m)})\!\rangle$ is a
constant~$A_m(\lambda)$.  These constants form a non-negative sequence defined
inductively by $A_0(\lambda) = 0$ and
$$
A_{m+1}(\lambda) = \frac{b_m(\lambda) (m+n-2)(n+2m)}{(m+1)} A_m(\lambda).
$$

Of course, since $K$ is positive, the sequence~$b_m(\lambda)$ must eventually
become negative.  If $\lambda$ were not of the form $m(n{+}m{-}1)K$, then this
sequence would become negative without ever becoming zero, which would force
the sequence~$A_m(\lambda)$ to become negative at some point.  This is
manifestly impossible, so there must exist an~$m$ so that
$\lambda=m(n{+}m{-}1)K$.  Of course, this implies
$$
0 = A_{m+1}(\lambda) = \langle\!(f^{(m+1)},f^{(m+1)})\!\rangle,
$$
which in turn implies that $f^{(m+1)}\equiv 0$.  By the Lemma, this implies
that $f$ is the restriction to $\Sigma$ of a harmonic polynomial of degree~$m$,
so that $f$ extends to the entire $n$-sphere, as desired.  This completes the
proof of Theorem~3.

\subhead 3.6. Non-uniqueness\endsubhead
For simplicity, I will now restrict to the case $K=1$ and assume that $\Sigma^n
= S^n$.  In~\cite{Br}, a strong uniqueness theorem was proved:  If $f:S^2\to
S^N$ is harmonic and of constant energy density, then $f$ is actually a minimal
isometric immersion up to a constant scale factor, thus implying that
$f\bigl(S^2)$ is a Boruvka sphere.  For $n$ greater than~2, and $N$ and $m$
sufficiently large, this sort of uniqueness fails, as can already be seen
in~\cite{DW} and ~\cite{DZ}.
In fact, the non-uniqueness in the present case of harmonic mappings follows
from the same argument: As has been shown, $f:S^n\to S^N$ is harmonic of
constant energy density~$\lambda = m(m{+}n{-}1)$ if and only if $f$ is the
restriction to~$S^n$ of a mapping $F:\bbR^{n+1}\to \bbR^{N+1}$ whose components
are harmonic polynomials of degree~$m$.  In fact, choosing a sequence
$F^0,\ldots, F^N$ of functions in~$\cH_m\bigl(\bbR^{n+1}\bigr)$, they are the
components of such an~$F$ if and only if
$$
(F^0)^2 + (F^1)^2 +\cdots + (F^N)^2 = \bigl((x^0)^2 + \cdots +
(x^n)^2\bigr)^{m} = R^m
$$
Let $D = \dim \cH_m\bigl(\bbR^{n+1}\bigr)$ and let $\{h_\alpha\,\vrule\,1\le
\alpha\le D\}$ be an orthonormal basis for $\cH_m\bigl(\bbR^{n+1}\bigr)$ with
respect to the inner product defined in~\S3.2.  Given any symmetric $D$-by-$D$
matrix $G= (g^{\alpha\beta})$, the expression $h(G) = g^{\alpha\beta}\,
h_\alpha h_\beta$ is a homogeneous polynomial of degree~$2m$ on~$\bbR^{n+1}$.
Note that, by equivariance, the identity matrix $G_0 = (\delta^{\alpha\beta})$
yields $h(G_0) = c R^m$ where $c>0$ is some constant.  Since, for large
enough~$m$, one has
$$
\dim S^d(\bbR^{n+1}) = {{n+2d}\choose{n}}
\ll \dim S^2\bigl(\cH_d(\bbR^{n+1})\bigr)
\approx \frac{2 d^{2n-2}}{\bigl((n-1)!\bigr)^2},
$$
(this strong inequality is where $n>2$ is used) it follows that the linear
mapping $G\mapsto h(G)$ cannot be injective.  The inverse image $\cG =
h^{-1}\bigl(R^m\bigr)$, which is a non-empty affine subspace since it
contains~$c^{-1}G_0$, is therefore of positive dimension and contains a closed
convex set~$\cC$ with non-empty interior consisting of positive semi-definite
matrices.

Each matrix $G\in\cG$ gives rise to a harmonic mapping of constant energy
density as follows.  Write $G = (g^{\alpha\beta}) = s^2$  where $s$ is also
positive semi-definite matrix (of the same rank as~$G$.  Then the functions
$F^\alpha = s^{\alpha\beta} h_\beta$ satisfy the above quadratic relation and
furnish a harmonic map $f:S^n\to S^{r-1}$ where $S^{r-1}\subset S^{D-1}$ is the
intersection of the unit sphere in~$\bbR^D$ with the $r$-dimensional range
of~$G$.

Again, since $n>2$, the dimension of~$\cC$ grows rather rapidly with~$m$, so
that it quickly exceeds the dimension of~$\SO(n{+}1)$, the group of domain
symmetries of the problem. At that point, non-uniqueness of $f$ up to the
obvious symmetries of the problem is assured.

\subhead 3.7. Non-positive sectional curvature\endsubhead
In \cite{Br}, where the $n=2$ case of Theorem~3 was treated, formulae like~(1)
were used in the case $K=0$ to prove that any harmonic mapping of constant
energy density from an open set in the Euclidean plane to an $N$-sphere was the
restriction of a harmonic mapping defined on the entire plane and was, in fact,
given by trigonometric polynomials.  It would be interesting to know if this
result holds in higher dimensions.

In the case $K<0$, it was shown that there were no harmonic maps of constant
energy density from $\Sigma$ to $S^N$ for any~$N$, a result which is also
unknown for higher dimensions.  This result depended in a crucial way on the
growth of the sequence~$A_m(\lambda)$ and other sequences which are closely
related.  It is still possible that the result for higher dimensions will
follow from such growth estimates, but I do not see how to carry them out.


\Refs

\widestnumber\key{BCG}

\ref\key BCG \by R. Bryant, et al
\book Exterior Differential Systems
\bookinfo
\publ Springer-Verlag
\publaddr New York
\yr 1991
\endref

\ref\key Br \by R. Bryant
\paper Minimal surfaces of constant curvature in the $n$-sphere
\jour Transactions of the American Mathematical Society
\vol 290
\yr 1985
\pages 259--271
\endref

\ref\key Ca1 \by E. Calabi
\paper Minimal immersions of surfaces in Euclidean spheres
\jour Journal of Differential Geometry
\vol 1
\yr 1967
\pages 111--125
\endref

\ref\key Ca2 \by E. Calabi
\paper Extremal isosystolic metrics for compact surfaces
\miscnote preprint
\endref

\ref\key DW \by M. do Carmo and N. Wallach
\paper Minimal immersions of spheres into spheres
\jour Annals of Mathematics
\vol 93
\yr 1971
\pages 43--62
\endref

\ref\key DZ \by D. DeTurck and W. Ziller
\paper Spherical minimal immersions of spherical space forms
\jour Comment. Math. Helv.
\vol 67
\yr 1992
\pages 428--458
\endref


\ref\key Gr \by M. Gromov
\paper Filling Riemannian manifolds
\jour Journal of Differential Geometry
\vol 18
\yr 1983
\pages 1--147
\endref

\ref\key Ka \by H. Karcher
\paper Embedded minimal surfaces derived from Scherk's examples
\jour Manuscr. Math.
\vol 63
\yr 1988
\pages 83-114
\endref

\ref\key Ke \by K. Kenmotsu
\paper Minimal surfaces of constant curvature in 4-dimensional space forms
\jour Proc. Amer. Math. Soc.
\vol 89
\yr 1983
\pages 133--138
\endref

\ref\key L \by P. Li
\paper Minimal immersions of compact irreducible homogeneous Riemannian
manifolds
\jour Journal of Differential Geometry
\vol 16
\yr 1981
\pages 105--115
\endref

\ref\key PU \by J.-S. Park and H. Urakawa
\paper Classification of harmonic mappings of constant energy density into
spheres
\jour Geometriae Dedicata
\vol 37
\yr 1991
\pages 211--226
\endref


\ref\key T \by T. Takahashi
\paper Minimal immersions of Riemannian manifolds
\jour J. Math. Soc. Japan
\vol 18
\yr 1966
\pages 380--385
\endref

\ref\key W \by N. Wallach
\paper Extension of locally defined minimal immersions of spheres into spheres
\jour Arch. Math.
\vol 21
\yr 1970
\pages 210--213
\endref

\endRefs
\enddocument